\documentclass[11pt,reqno]{amsart}
\usepackage[body={7in,9in},centering]{geometry}
\usepackage{amsmath,amssymb,amsthm,mathrsfs}
\usepackage{color,xcolor}
\usepackage{setspace}
\usepackage{parskip}
\definecolor{cobalt}{RGB}{61,99,181}

\usepackage[colorlinks,linkcolor=blue,anchorcolor=blue,citecolor=blue]{hyperref}

\newtheorem{thm}{Theorem}[section]
\newtheorem{cor}[thm]{Corollary}

\newtheorem{ex}[thm]{Example}

\numberwithin{equation}{section}

\date{\today}

\makeatletter

\newcommand{\Rmnum}[1]{\expandafter\@slowromancap\romannumeral #1@}
\makeatother

\begin{document}

\title[DYNAMIC OF TIME VARYING WEIGHTED BACKWARD SHIFT OPERATOR  ]{DYNAMIC OF TIME VARYING WEIGHTED BACKWARD SHIFT OPERATOR }

\author[Jing Hou]{Jing Hou\textsuperscript{1}}
\address{\textsuperscript{1} College of Mathematics and Statistics, Chongqing University, Chongqing, 401331, P. R. China}
\email{Jing Hou@cqu.edu.cn}
\author[yonglu Shu]{yonglu Shu\textsuperscript{2}}
\address{\textsuperscript{2} College of Mathematics and Statistics, Chongqing University, Chongqing, 401331, P. R. China}
\email{shuyonglu@163.com}


\begin{abstract}
we study the hypercyclic and chaotic properties of the time varying weighted backward shift operator $(Tx)(t)=w(t)x(t+a)$ in $L_p(0,\infty)(1\leq p<\infty)$ and $C_0[0,\infty)$. And we also analyse the spectral structure of the operators if the spaces are complex.
\end{abstract} \maketitle
\section{INTRODUCTION}\label{S1}
Linear dynamics is a new branch of functional analysis which has developed rapidly. The study of linear dynamics began in 1929 when Birkhoff gave the hypercycle of translation operator in $H(\mathbb{C})$ \cite{birkhoff1929demonstration} . Then, in 1952, Maclane discovered that differential operator in $H(\mathbb{C})$ is hypercyclic \cite{maclane1952sequences}. In 1969, Rolewicz showed that linear shifts can have dense orbits in Banach space \cite{rolewicz1969orbits}. However, the three papers mentioned above have not attracted much attention for mathematicians.
 The real development of linear operator dynamics was after the publication of Kaitai's thesis \cite{kitai1984invariant} in 1984, which discussed the problem of invariant subsets of linear operators and gave the Hypercyclic criterion. Since the 1980s, many mathematicians, inspired by the above sporadic examples, have aroused the enthusiasm of researchers and begun to study the dynamical properties of general linear operators, see \cite{godefroy1991operators}, \cite{grosse1999universal}, \cite{shapiro2001notes}.

We refer to the books \cite{bayart2009dynamics} or \cite{LinearChaos} for background on linear dynamics, and to the papers \cite{glasner2015universal} or \cite{grivaux2017some} for a glimpse at the richness of the class of linear dynamical systems and its potential usefulness in general ergodic theory. If T is a bounded operator, the $T-$orbit of a vector $x\in X$ is the set $orb(x,T):=\{T^nx\}_{n\in \mathbb{Z}_+}$. The operator $T$ is said to be $hypercyclic$ if there is some vector $x\in X$ such that $orb(x,T)$ is dense in $X$. Such a vector $x$ is said to be hypercyclic for $T$, and the set of all hypercyclic vectors for $T$ is denoted by $HC(T)$. $Hypercyclic$ linear operators with a dense set of periodic points are said to be $chaotic$. In paper \cite{bes2001chaotic}, let $A$ be an unbounded operator on a separable infinite dimensional Banach space $X$. Since $A$ is only densely defined on $X$, it may happen that a vector $f$ is in the domain of $A$, but $Af$ fails to be in the domain of $A$. For this reason, we define $A$ to be $hypercyclic$ if there is a vector $f$ in the domain of $A$ such that for every integer $n>1$ the vector $A^nf$ is in the domain of $A$ and the orbit $\{f, Af, A^2f, A^3f,\cdots\}$ is dense in $X$. $Hypercyclic$ unbounded operators with a dense set of periodic points are said to be $chaotic$.

In \cite{rolewicz1969orbits}, S. Rolewicz provides the first example of a hypercyclic bounded linear operator on a Banach space $l_p (1\leq p<\infty)$ or $c_0$ weighted backward shift:$T(x_k)_{k\in\mathbb{N}}=w(x_{k+1})_{k\in\mathbb{N}}$, where $w\in\mathbb{F}(\mathbb{F}=\mathbb{R}~\text{or}~\mathbb{C})$ such that $|w|>1$ is arbitrary. Rolewicz's operators also happen to chaotic.

In 2019, Marat V. Markin considers a natural unbounded extension of the classcial Rolewicz weighted backward shift operators\cite{markin2019real}. He prove the chaoticity and describe the spectral structure of Rolewicz-type weighted backward shift unbounded linear operator in the sequence spaces $l_p(1\leq p<\infty)$ and $c_0$: $A(x_k)_{k\in\mathbb{N}}=(w^kx_{k+1})_{k\in\mathbb{N}}$, where $w\in\mathbb{R}$ or $w\in\mathbb{C}$ such that $|w|>1$ is arbitrary, with the domain $D(A)=\{(x_k)_{k\in\mathbb{N}}\in X| (w^kx_k+1)_{k\in\mathbb{N}}\in X\}$ where $X=l_p(1\leq p<\infty)$ or $c_0$.

In \cite{jimenez2022linear}, John M. Jimenez and Marat V. Markinn shown that, in $L_p(0,\infty)$ and $C_0[0,\infty)$, the bounded weighted backward shift operator $(Tx)(t)=wx(t+a)~(w>1~\text{and}~a>0)$ and its unbounded counterpart $(Tx)(t)=w^tx(t+a)$ are chaotic. Combining \cite{jimenez2022linear} and \cite{markin2018chaoticity}, a natural question is : what are the dynamical properties if we replacing $w^t$ and $w$ with $w(t)$?

 We can define shift operators $T$ with variable weight $w(t)$:~$Tx(t)=w(t)x(t+a)$.
We get the main theorems of this paper:

\begin{thm}\label{1}In the $($real or complex$)$ space $L_p(0,\infty) (1\leq p<\infty)$ the bounded time varying weight backward shift operator
\begin{align*}
Tx(t):&=w(t)x(t+a),\\
T^nx(t):&=\prod_{i=0}^{n-1}w(t+ia)x(t+na),
\end{align*}
where w(t) bounded and $|w(t)|>1$, is chaotic. Moreover, supposed the space is complex,$$\sigma(T)=\{\lambda\in\mathbb{C}|~|\lambda|\leqslant\|w\|_\infty\}$$
with $$\sigma_p(T)=\{\lambda\in\mathbb{C}|~|\lambda|<\|w\|_\infty\}~\text{and}~\sigma_c(T)=\{\lambda\in\mathbb{C}|~|\lambda|=\|w\|_\infty\}.$$
\end{thm}

\begin{thm}\label{2}
In the $($real or complex$)$ space $L_p(0, \infty)$, we can define the unbounded time varying weight backward shift operator $T$
\begin{align*}
(Tx)(t)&:=w(t)x(t+a)\\
(T^nx)(t)&:=\prod\limits_{i=0}^{n-1}w(t+ia)x(t+na),
\end{align*}
where $x\in L_p(0, \infty)$ and $a>0$, with domain $D(T):=\{x\in L_p(0, \infty)|\int_0^\infty|w(t)x(t+a)|^pdt<+\infty\}.$
If $w(t)$ satisfies the following conditions$:$

$(1)w(t)$ unbounded and

$(2)\lim\limits_{t\rightarrow\infty}w(t)=\infty$,

$(3)|w(t)|>b>1,$ $($b is a constant$)$.\\
then $T$ is chaotic. Moreover, if the space is complex, then$$\sigma(T)=\sigma_p(T)=\mathbb{C}.$$

\end{thm}

\begin{thm}\label{3}
In the $($real or complex$)$ space $C_0[0, \infty)$, we can define the unbounded time varying weighted backward shift operator $T$
\begin{align*}
(Tx)(t)&:=w(t)x(t+a)\\
(T^nx)(t)&:=\prod\limits_{i=0}^{n-1}w(t+ia)x(t+na),
\end{align*}
where $x\in C_0[0, \infty)$ and $a>0$, with domain $$D(T):=\{x\in C_0[0, \infty)|~\lim_{t\rightarrow\infty}|w(t)x(t+a)=0\}.$$
If $w(t)$ satisfies the following conditions$:$

$(1)$ w(t) is continuous, $w(0)=1$ and $\lim\limits_{t\rightarrow\infty} w(t)=\infty$;

$(2)|w(t)|>b>1,$ $($b is a constant$)$ and

$(3)\max\limits_{t\in[0,\infty)}|\frac{w(t)}{w(t-(n-1)a)}|\rightarrow\infty,~ n\rightarrow \infty,$\\
then $T$ is hypercyclic.
\end{thm}

\section{PRELIMINARY}\label{S2}
In this section, we recall some well-known definitions, referring to \cite{bayart2009dynamics} or \cite{LinearChaos} for more retails. And all operators are considered to be linear operators in a separable Banach space $X$ with
norm $\|\cdot\|$.
\begin{cor}(Hypercyclicity\cite{bayart2009dynamics}\cite{LinearChaos}).
Let $$T:X\rightarrow X$$
be a bounded operator in a Banach space$(X,\|\cdot\|)$. A bounded operator $T$ possessing hypercyclic vectors is called $hypercyclic$.

 A nonzero vector $x\in X$ is called $hypercyclic$ for $T$ if its orbit under $T$, $$orb(x,T)=\{T^nx;~n\geqslant0\}$$ is dense in $X$.
\end{cor}
\begin{cor}(chaotic\cite{bayart2009dynamics}\cite{LinearChaos}). A periodic point for $T$ is a vector $x\in X$ such that $T^nx=x$ for some $n\in \mathbb{N}$.
Provided that a bounded operator $T$ satisfies the following conditions:\\
(i) $T$ is hypercyclic;\\
(ii)$T$ has a dense set of periodic points,\\
then $T$ is chaotic.
\end{cor}
\begin{cor} For an operator $T$ in $X$ with domain $D(T)$, the set $D(T^\infty)$ is by definition the intersection of all $D(T^n)$'s: $$D(T^\infty):=\bigcap_{n=0}^\infty D(T^n).$$
with $\|x\|_n:=\sum_{j=0}^{n}\|T^jx\|$ for $x\in D(T^\infty)$.
\end{cor}
In \cite{delaubenfels2003chaos}, if $T$ is a closed operator then the space $D(T^\infty)$ is a separable Frechet space under the topology determined by the seminorms $\|x\|_n, n=0, 1,  2, \cdots$. There is a natural extension of the concept of hypercyclicity and chaos to unbounded operators, see \cite{bes2001chaotic} and \cite{delaubenfels2003chaos}.

\begin{cor}(Hypercyclicity \cite{bes2001chaotic}\cite{delaubenfels2003chaos}).
Let $$T:D(T)\rightarrow X$$
be a unbounded operator in a Banach space$(X,\|\cdot\|)$. A unbounded operator $T$ possessing hypercyclic vectors is called $hypercyclic$.

 A nonzero vector $x\in D(T^\infty)$ is called $hypercyclic$ for $T$ if its orbit under $T$, $$orb(x,T)=\{T^nx;~n\geqslant0\}$$ is dense in $X$.
\end{cor}
\begin{cor}(chaotic\cite{bes2001chaotic}\cite{delaubenfels2003chaos}). A periodic point for $T$ is a vector $x\in D(T^\infty)$ such that $T^nx=x$ for some $n\in \mathbb{N}$.
Provided that an operator $T$ satisfies the following conditions:\\
(i) $T$ is hypercyclic;\\
(ii)$T$ has a dense set of periodic points,\\
then $T$ is chaotic.
\end{cor}
For bounded linear operators, theorem \ref{chuandidingli} is a convenient tool for establishing their hypercyclicity.
\begin{thm}\label{chuandidingli}(Birkhoff transitivity theorem\cite{LinearChaos}). Supposing a bounded linear operator $T$ on a separable Banach space $(X,\|\cdot\|)$ , we set, for any subsets $A,~B$ of $X$,$$N_T(A,B):=\{i\in\mathbb{N}|~T^i(A)\cap B\neq\emptyset\}.$$ If a bounded operator $T$ is hypercyclic if and only if $T$ is topologically transitive, i.e. $N_T(U,V)\neq\emptyset$ for all open sets $U,~V\neq\emptyset$.
\end{thm}

  We find that we use the above theorem for bounded linear operators, but for unbounded linear operators, we use

\begin{thm}\label{chaoxunhuan}(Sufficient Condition for Hypercyclicity\cite{bes2001chaotic})
Let $(X,\|\cdot\|)$ be a separable infinite-dimensional Banach space and let $(T,D(T))$ be a densely defined linear operator in $X$ for which the operator $T^n$ is closed for all $n\in\mathbb{N}$.\\
Then the operator $T$ is called hypercyclic if there exists $$\epsilon\subset D(T^\infty):=\bigcap_{n=1}^\infty D(T^n) $$ dense in $(X,\|\cdot\|)$ and a mapping $ S:\epsilon\rightarrow\epsilon$ such that \\

(i) $TS=I$, (I is the identity mapping on $\epsilon$) and

(ii)$\forall x\in\epsilon,~T^nx,S^nx\rightarrow0,~n\rightarrow\infty$, in $(X, \|\cdot\|)$.
\end{thm}

\begin{cor}(Spectrum Classification\cite{dunford1988linear}\cite{markin2020elementary}). The spectrum $\sigma(T)$ of a closed linear operator (T,D(T))in a complex Banach space $(X,\|\cdot\|)$ is the union of the following three pairwise disjoint sets:

$\sigma_p(T)=\{\lambda\in\mathbb{C}|~T-\lambda I\text{is not injective}\}$,

$\sigma_c(T)=\{\lambda\in\mathbb{C}|~T-\lambda I\text{is injective and not surjective, but }\overline {R(T-\lambda I)}=X\},$

$\sigma_r(T)=\{\lambda\in\mathbb{C}|~T-\lambda I\text{is injective and }\overline{R(T-\lambda I)}\neq X\}$,\\
called the point, continuous, and residual spectrum of $T$.
\end{cor}
\section{ Proof of Theorem \ref{1}}
In this section, we will prove the Theorem \ref{1}. The tool used to prove that the bounded time varying weighted backward shift operator is  hypercycle is Theorem \ref{chuandidingli}.
\begin{proof}
Let~$U, V\subset L_p(0, \infty)$~be arbitrary nonempty open set. By the denseness in ~$L_p(0, \infty)$~of $p$ integrable eventually zero function, there exist~$x\in U$~and~$y\in V$~. Since these functions are eventually zero, we can find a sufficiently large $N\in\mathbb{N}$ such that $$x(t)=0~\text{and}~y(t)=0,~t\geqslant Na.$$
For an arbitrary $n\geqslant N$, the $p$ integrable eventually zero function
$$z_n(t):=\begin{cases}
x(t) &t\in[0, Na),\\
\big[\prod\limits_{i=0}^{n-1}w(t-ia+na)\big]^{-1}y(t-an)&t\in[na, Na+na),\\
0&\text{otherwise,}
\end{cases}$$
~in~$L_p(0,\infty).$

Observe, if $t\in[0, Na)$ \begin{align*}T^nz_n(t)&=\prod_{i=0}^{n-1}w(t+ia)z_n(t+na)\\
&=\prod_{i=0}^{n-1}w(t+ia)[\prod_{i=0}^{n-1}w(t-ia)]^{-1}y(t)\\
&=y(t).
\end{align*}
If $t\geq Na$, \begin{align*}
T^nz_n(t)&=\prod_{i=0}^{n-1}w(t+ia)z_n(t+na)\\
&=0\\
&=y(t).
\end{align*}
Then
$$T^nz_n(t)=y(t), t\geq 0.$$
Since $w(t)$ is a bounded function, then $\exists~m>1$ such that $|w(t)|\geq m$.
Then \begin{align*}\|z_n-x\|_p&=\big(\int_0^{\infty}|z_n(t)-x(t)|^pdt\big)^{\frac{1}{p}}\\
&=\big(\int_{na}^{Na+na}\big{|}\frac{1}{\prod_{i=0}^{n-1}w(t+ia-na)}y(t-na)\big |^pdt\big)^{\frac{1}{p}}\\
&\leqslant \frac{1}{m^n}\|y\|_p\rightarrow 0, n\rightarrow\infty.
\end{align*}

Consequently, for all sufficiently large $n\in\mathbb{N}$,~$z_n\in U$~and~$T^nz_n=y\in V$. Thus $T$ is topologically transitive.
By the Birkhoff Transitivity Theorem, $T$ is hypercyclic.

Now, we are going to show that  T has a dense set of periods points. By the denseness in $L_p(0, \infty)$ of the equivalence classes represented by $p$-integrable eventually zero functions, it suffices to approximate any such a class by periodic points of T.

For any $x\in L_p(0, Na)$ represented by a $p$-integrable on $(0, Na)$ function $x(\cdot)$ and $\forall N\in \mathbb{N}$, the following $p-$integrable on $(0, \infty)$ function
\begin{equation}x_N(t):=\frac{1}{\prod\limits_{i=1}^{kN}w(t-ia)}x(t-kNa), t\in D_k,\end{equation}
where $D_k=[kNa,(k+1)Na), k\in Z_+$. Since
 \begin{align*}T^Nx_N(t)&=\prod\limits_{i=0}^{N-1}w(t+ia)x_N(t+Na)\\
 &=\prod\limits_{i=0}^{N-1}w(t+ia)\frac{1}{\prod\limits_{i=1}^{kN}w(t+Na-ia)}x(t-(k-1)Na)\\
 &=\frac{1}{\prod\limits_{i=1}^{(k-1)N}w(t-ia)}x(t-(k-1)Na).
\end{align*}
 Thus $T^Nx_N(t)=x_N(t)$, then we can get $x_N(T)$ is a periodic point of $T$ with period $N$.

 Now we can show that $x_N\in L_p(0, \infty)$. Indeed,
 \begin{align*}\int_0^\infty|x_N(t)|^pdt&=\int_0^\infty|[{\prod\limits_{i=1}^{kN}w(t-ia)}]^{-1}x(t-kNa)|^pdt\\
&=\sum\limits_{k=0}^\infty\int_{D_k}\frac{1}{\prod\limits_{i=1}^{kN}|w(t-ia)|^p}|x(t-kNa)|^pdt\\
&\leqslant\sum\limits_{k=0}^\infty\frac{1}{m^{kNp}}\int_{D_k}|x(t-kNa)|^Pdt\\
&=\sum\limits_{k=0}^{\infty}\frac{1}{m^{kNp}}\int_{0}^{Na}|x(t)|^pdt.
 \end{align*}
Since $m>1$ and $x(t)\in L_p(0, Na)$, then $\int_{0}^{\infty}|x(t)|^pdt<\infty$. Thus $x_N(t)\in L_p(0, \infty).$

Now, let $\forall~y\in L_p(0, \infty)$, then $\exists~n\in\mathbb{N},~y(t)=0,~t>na.$~Let $y_N$ be the periodic point of $T$ of an arbitrary period $N>n$, defied base on $y(\cdot)$ as in (3.1).

Thus,
\begin{align*}
\|y-y_N\|_p&=(\int_0^\infty|y(t)-y_N(t)|^pdt)^{\frac{1}{p}}\\
&=\sum\limits_{k=0}^\infty(\int_{D_k}|y(t)-y_N(t)|^pdt)^{\frac{1}{p}}\\
&=(\int_0^{Na}|y-y_N(t)|^pdt)^{\frac{1}{p}}+\sum\limits_{k=1}^{\infty}(\int_{D_k}|y(t)-y_N(t)|^{p})^{\frac{1}{p}}\\
&=0+\sum\limits_{k=1}^{\infty}(\int_{D_k}|y_N(t)|^{p})^{\frac{1}{p}}\\
&\leqslant\sum\limits_{k=1}^{\infty}m^{-kN}\|y\|_p\\
&=\frac{m^{-N}}{1-m^{-N}}\|y\|_p\rightarrow0,~N\rightarrow\infty.
\end{align*}

This implies that the set $Per(T)$ of periodic points of T in dense in $L_p(0, \infty)$, and hence the operator is chaotic.

Now, suppose that the space $L_p(0,\infty) $ is complex, we analyze the spectrum of operator $T$.

Firstly, we show that $\|T\|=\|w\|_\infty$.
We can get that $\|T\|\leqslant\|w\|_\infty$, since\begin{align*}\|T\|=\sup_{\|x\|_p=1}\|Tx\|_p&=\|w\|_\infty\sup_{\|x\|_p=1}\bigg[\int_0^{+^\infty}|x(t+a)|^p~dt\bigg]^{\frac{1}{p}}\\
                                                  &\leqslant\|w\|_\infty\sup_{\|x\|_{p=1}}\bigg[\int_0^{+\infty}|x(t)|^p~dt\bigg]^{\frac{1}{p}}\\
                                                  &\leqslant\|w\|_\infty\sup_{\|x\|_{p=1}}\bigg[\int_0^{+\infty}|x(t)|^p~dt\bigg]^{\frac{1}{p}}\\
                                                  &=\|w\|_\infty.
     \end{align*}
We can find that operator $T$ can decomposed into the product of  $V$ and $\tau_a$, where $V: x(t)\mapsto w(t)x(t)$ and $\tau_a: x(t)\mapsto x(t+a)$. For $x\in L_p(0,+\infty)$, we have $$T\tau_{-a}x(t)=Tx(t-a)=w(t)x(t)=Vx(t),$$ then $T\cdot\tau_{-a}=V$. Thus $\|V\|=\|V\|\leqslant\|T\|\|\tau_{-a}\|=\|T\|$ where $\|\tau_{-a}\|=1.$ Then we can get $\|T\|=\|w\|_{\infty}.$ By Gelfand's Spectral Radius Theorem for bounded linear operator, we have $$\sigma(T)\subseteqq\{\lambda\in\mathbb{C}|~|\lambda|\leqslant\|w\|_\infty\}.$$

Secondly, we show that $$\sigma_p(T)=\{\lambda\in\mathbb{C}|~|\lambda|<\|w\|_\infty\}.$$

(i) Let $\lambda\in\mathbb{C}$, with $|\lambda|<\|W\|$ and a nonzero equivalence class $\forall y \in L_p(0,a)$. Consider the function $$x(t):=\frac{\lambda^k}{w(t-a)w(t-2a)\cdots w(t-ka)}y(t-ka),$$where $t\in[ka,(k+1)a),~k\in\mathbb{Z}_+.$

Then $x\in D(T)$ and $Tx=\lambda x.$ Indeed, since $$\frac{y(t-ka)}{w(t-a)w(t-2a)\cdots w(t-(k-1)a)}\in L_p(0,\infty)\setminus\{0\},$$
 then we can get  $$x(t)\in L_p(0,\infty)\setminus\{0\}.$$
Also \begin{align*}
(Tx)(t)&=w(t)x(t+a)\\
       &=w(t)\frac{\lambda^k}{w(t)w(t-a)\cdots w(t-(k-1)a)}y(t-(k-1)a)\\
       &=\lambda\frac{\lambda^{k-1}}{w(t-a)w(t-2a)\cdots w(t-(k-1)a)}y(t-(k-1)a)\\
       &=\lambda x(t),
\end{align*}
which implies that $$Tx=\lambda x.$$
This shows that $$\{\lambda\in\mathbb{C}|~|\lambda|<\|w\|_\infty\}\subset\sigma_p(T).$$

(ii) Conversely, we show that $$\sigma_p(T)\subset\{\lambda\in\mathbb{C}|~|\lambda<\|w\|_\infty|\}.$$
For $x\in D(T)\backslash\{0\}$ where $\forall\lambda\in\mathbb{C}$ satisfies $Tx=\lambda x.$

Let $$x_k(t):=x(t),~t\in[ka,(k+1)a),~k\in \mathbb{Z}_+.$$
We get that $$\lambda x_{k-1}(t)=w(t)x_k(t+a),~t\in[ka,(k+1)a),~k\in\mathbb{Z}_+.$$
Hence$$x_k(t)=\frac{\lambda^{k}}{w(t-a)w(t-2a)\cdots w(t-ka)}y(t-ka),~t\in[ka,(k+1)a),~k\in\mathbb{Z}_+,$$
which implies that$$0<\int_{0}^{a}|x_k(t)|^p~dt\leqslant\int_0^\infty|x_k(t)|^p~dt<+\infty$$
and \begin{align*}
\infty>\int_{0}^{\infty}|x(t)|^p~dt&=\sum_{k=0}^{\infty}\int_{ka}^{(k+1)a}\bigg|\frac{\lambda^k}{w(t-a)w(t-2a)\cdots w(t-ka)}y(t-ka)\bigg|^p~dt\\
&\geqslant \sum_{k=0}^{\infty}\bigg(\frac{|\lambda|}{\|w\|_{\infty}}\bigg)^{kp}\int_{ka}^{(k+1)a}|y(t-ka)|^p~dt\\
&=\sum_{k=0}^{\infty}\bigg(\frac{|\lambda|}{\|w\|_\infty}\bigg)^{kp}\int_{0}^{a}|y(t)|^p~dt.
\end{align*}
Since $x\in L_p(0,a)$, we must satisfy $$\sum\limits_{k=0}^{\infty}\bigg(\frac{|\lambda|}{\|w\|\small_\infty}\bigg)^{kp}<+\infty$$ for $$\sum_{k=0}^{\infty}\bigg(\frac{|\lambda|}{\|w\|_\infty}\bigg)^{kp}\int_{0}^{a}|y(t)|^p~dt<+\infty.$$
Hence we can imply that $$\lambda<\|w\|_\infty$$ and $$\sigma_p(T)\subset\{\lambda\in\mathbb{C}|~|\lambda|<\|w\|_\infty\}.$$
Considering that $\sigma(T)$ is a closed subset in $\mathbb{C}$, we get that $$\sigma(T)=\{\lambda\in\mathbb{C}|~|\lambda|\leqslant\|w\|_\infty\}.$$
Since $T-\lambda I$ has dense range, which shows that $$\sigma_r(T)=\emptyset.$$
It easily to show that $$\sigma_c(T)=\{\lambda\in\mathbb{C}|~|\lambda|=\|w\|_\infty\},$$ by $\sigma(T)=\sigma_p(T)\cup\sigma_c(T)\cup\sigma_r(T)$ where $\sigma_p(T),\sigma_c(T),\sigma_r(T)$ pairwise disjoint. This completes the proof.
\end{proof}
\section{Proof of Theorem \ref{2}}\label{S4}
In Theorem \ref{2}, the operator is unbounded. Then, if we prove that the unbounded time varying weight shift operator is hypercycle, Theorem \ref{chuandidingli} will not be used, but Theorem \ref{chaoxunhuan} will be used. Now let's prove the theorem \ref{2}.
\begin{proof}
Firstly, we show that T is unbounded.\\Consider $$x_n(t):=\chi_{[n, n+1]}(t)=\begin{cases}1,~~t\in[n, n+1],\\0,~~otherwise.\end{cases}$$
$x_n\in D(T)$ and $\|x_n\|_p=1, n\in \mathbb{N}.$

Then, for all sufficiently large $n\in \mathbb{N}$,
\begin{align*}
\int_0^\infty|Tx_n(t)|^pdt&=\int_0^\infty|w(t)x_n(t+a)|^pdt\\
&=\int_{n-a}^{n+1-a}|w(t)|^pdt\rightarrow\infty,~n\rightarrow\infty.
\end{align*}
Which implies that $T$ is unbounded.

Secondly, we will show that $T^l$ is closed for each $l\in \mathbb{N}.$ Let a sequence $(x_n)_{n\in \mathbb{N}}$ be such that
$$D(T)\ni x_n\rightarrow x\in L_p(0, \infty),~n\rightarrow\infty$$
and
$$T^lx_n\rightarrow y\in L_p(0, \infty),~n\rightarrow\infty.$$

Then, the function sequences $(x_n(\cdot))_{n\in\mathbb{N}}$ and $(T^lx_n(\cdot))_{n\in\mathbb{N}}$ being convergent in p-norm, converge in the Lebesgue measure on $(0,\infty)$,
and hence, by the Riesz theorem, there exist subsequences $(x_{n_k}(\cdot))_{k\in\mathbb{N}}$ and $(T^lx_{n_k}(\cdot))_{k\in\mathbb{N}}$ such that
\begin{align*}
x_{n_k}(t)&\rightarrow x(t),~k\rightarrow\infty,~a.e. t\in(0,\infty),\\
T^lx_{n_k}(t)&\rightarrow y(t),~k\rightarrow\infty,~a.e. t\in(0,\infty).
\end{align*}
Thus, $$\prod_{i=0}^{l-1}w(t+ia)x_{n_k}(t+la)\rightarrow\prod_{i=0}^{l-1}w(t+ia)x(t+la),~a.e.~t\in(0,\infty).$$
By the $completeness$ of the $Lebesgue$ measure, we infer that
Which implies that $x\in D(T^l)~\text{and}~T^lx=y$. It following that $T^l$ is closed. Next, we are going to show that $T$ is $hypercyclic.$

Let $\epsilon$ be the set of the equivalence classes represented by p-integrable eventually zero functions. Note that $$\epsilon\subseteq C^{\infty}(T).$$
Consider $$\epsilon\ni x\mapsto(Sx)\in\varepsilon$$
where the $Sx$ is represented by$$(Sx)(t):=\begin{cases}\frac{1}{w(t-a)}x(t-a),~&t>a,\\0,~&otherwise.\end{cases}$$
Since $Sx(\cdot)$ is eventually zero and
\begin{align*}
\int_{0}^{\infty}|Sx(t)|^p~dt&=\int_{a}^\infty|Sx(t)|^p~dt\\
                             &=\int_{a}^{\infty}\bigg|\frac{1}{w(t-a)}x(t-a)\bigg|^p~dt\\
                             &=\int_{0}^{\infty}\bigg|\frac{1}{w(t)}x(t)\bigg|^p~dt\\
                             &\leqslant\int_{0}^{\infty}|x(t)|^p~dt.
\end{align*}
Thus the mapping $S:\epsilon\rightarrow\epsilon$ is well defined.

we can easy to show that $TS=I$. For any $x\in \varepsilon$, we know that $$supp~x=\{t\in(0,\infty)|x(t)\neq0\}$$
By $x\in C_0(0,\infty)$, we have that $$\exists~N\in\mathbb{N},~supp~x\subseteq[0,Na].$$
Because the definition of $T$, we can get that for any $n\geqslant N$, $$T^nx=0.$$
Therefore, it is proved that $\forall x\in\epsilon,~S^nx\rightarrow0,~n\rightarrow\infty.$

Next, we will show that for all $x\in\epsilon$, $$S^nx\rightarrow0,~n\rightarrow\infty.$$
When $\forall x\in\epsilon$, we can conclude that $\forall n\in\mathbb{N}$,
$$(S^nx)(t)=\begin{cases}\frac{1}{\prod\limits_{i=1}^{n}w(t-ia)}x(t-na),~&t>na,\\0,~&otherwise.
\end{cases}
$$
Because \begin{align*}
\|S^nx\|_p^p&=\int_{0}^\infty|(S^nx)(t)|^p~dt\\
            &=\int_{na}^{\infty}\bigg|\frac{1}{\prod\limits_{i=1}^{n}w(t-ia)}x(t-na)\bigg|^p~dt\\
            &\leqslant\int_{na}^\infty\bigg|\frac{1}{b^{np}}\bigg||x(t-na)|^p~dt\\
            &=\frac{1}{b^{np}}\int_0^{\infty}|x(t)|^p~dt\rightarrow0,~n\rightarrow\infty,
\end{align*}
we proved that  for all $x\in\epsilon$, $$S^nx\rightarrow0,~n\rightarrow\infty.$$

Therefore, by Theorem \ref{chaoxunhuan}, we get that $T$ is $hypercyclic$.

At present, we will prove that $T$ has a dense set of periodic point.
Since $C_0(0,\infty)$ is dense in $L_p(0,\infty)$, we only show that $per(T)$ is dense in $C_0(0,\infty)$.

Firstly, for $\forall x\in L_p(0,Na)$ and $\forall N\in\mathbb{N}$, we can define $x_N\in L_p(0,\infty)$ with the following representative
$$x_N(t):=\frac{1}{\prod\limits_{i=1}^{kN}w(t-ia)}x(t-kNa),~t\in[kNa,(k+1)Na),~k\in\mathbb{Z_+}.$$
It is easy to show that $x_N\in per(T)$.

Since, with $$D_k:=[kNa,(k+1)Na),~k\in\mathbb{Z_+},
$$
 \begin{align*}\int_0^\infty|x_N(t)|^pdt&=\int_0^\infty|[{\prod\limits_{i=1}^{kN}w(t-ia)}]^{-1}x(t-kNa)|^pdt\\
&=\sum\limits_{k=0}^\infty\int_{D_k}\frac{1}{\prod\limits_{i=1}^{kN}|w(t-ia)|^p}|x(t-kNa)|^pdt\\
&\leqslant\sum\limits_{k=0}^\infty\frac{1}{b^{kNp}}\int_{D_k}|x(t-kNa)|^Pdt\\
&=\sum\limits_{k=0}^{\infty}\frac{1}{b^{kNp}}\int_{0}^{Na}|x(t)|^pdt.
 \end{align*}
 Since $b>1$ and $x(t)\in L_p(0,Na)$, which implies that $x_N\in L_p(0,\infty).$ Therefore,$$\forall N\in\mathbb{N},~x_N\in L_p(0,\infty).$$
 If $y\in C_0(0,\infty)$, then $$\exists n\in\mathbb{N},~y(t)=0,~t>na.$$
 For an arbitrary period $N\geqslant n$, we defined $y_N$ with the following representative
$$y_N(t):=\frac{1}{\prod\limits_{i=1}^{kN}w(t-ia)}y(t-kNa),~t\in[kNa,(k+1)Na),~k\in\mathbb{Z_+}.$$
Now, we prove that $$\|y-y_N\|_p^p\rightarrow0,~N\rightarrow\infty.$$
Since$|w(t)|>b>1$, we have
\begin{align*}
\|y-y_N\|_p^p&=\int_0^\infty|y(t)-y_N(t)|^p~dt\\
             &=\sum_{k=1}^\infty\int_{D_k}\frac{1}{|\prod\limits_{i=1}^{kN}w(t-ia)|^p}|y(t-kNa)^p~dt\\
             &\leqslant\sum_{k=1}^\infty\int_{D_k}\frac{1}{b^{kNp}}|y(t-kNa)|^p~dt\\
             &=\sum_{k=1}^\infty\frac{1}{b^{kNp}}\int_{0}^{Na}|y(t)|^p~dt\rightarrow0,~N\rightarrow\infty.
\end{align*}
We have already shown that $per(T)$ is dense in $L_p(0,\infty)$. And, which implies that the operator $T$ is a chaotic.
Now, if the space $L_p(0,\infty)$ is complex, we analyse the spectrum of $T$.

We already know the fact $$\sigma_p(T)\subset\mathbb{C}.$$
Thus, we only show that $$\mathbb{C}\subset\sigma_p(T).$$
For $\forall\lambda\in\mathbb{C}$ and a nonzero $y\in L_p(0,a)$, we define that
$$x(t):=\frac{\lambda^n}{\prod\limits_{i=1}^nw(t-ia)}y(t-na),~t\in[na,(n+1)a),~n\in\mathbb{Z_+}.$$
Since\begin{align*}
\int_{0}^{\infty}|x(t)|^p~dt=\sum_{n=0}^{\infty}\int_{D_n}\bigg|\frac{\lambda^n}{\prod\limits_{i=1}^nw(t-ia)}y(t-na)\bigg|^p~dt<\infty.
\end{align*}
Which implies that $x\in L_p(0,\infty)\setminus\{0\}$ and
$$Tx=\lambda x.$$
Thus $\mathbb{C}\subset\sigma_p$, it is following that $$\sigma(T)=\sigma_p(T)=\mathbb{C}.$$
\end{proof}
Now we will give an example of unbounded time varying weight backward shift operator in $L_p(0,\infty), 1\leq p<\infty$.
\begin{ex}
In the (real or complex) space $L_p(0,\infty)$, we define that the unbounded weight backward shift operator T:
\begin{align*}
(Tx)(t)&:=[a_mt^m+a_{m-1}t^{m-1}+\cdots+a_0]x(t+a)\\
(T^nx)(t)&:=\prod_{i=0}^{n-1}P_m(t+ia)x(t+na),
\end{align*}
where $P_m(t)=a_mt^m+a_{m-1}t^{m-1}+\cdots+a_0,~m>1,~a_i>0,~i=0,1,2,\cdots,m$ and $a>0$, with domain $$D(T)=\bigg\{x\in L_p(0,\infty)\bigg|~\int_0^\infty|P_m(t)x(t+a)|^p~dt<\infty\bigg\},$$
is chaotic.
\end{ex}
\begin{proof}
We first show that T is unbounded.

Consider $$x_n(t):=\chi_{[n,n+1]}(t)=\begin{cases}1,~~~t\in[n,n+1],\\0,~~~otherwise\end{cases}$$as a representative of a class $x_n\in D(T)$. We can find that $$\|x_n\|_p=1,~ n\in\mathbb{N}.$$
Then, since $a_i>0,~i=1,2,\cdots,m,$ for all sufficiently large $n\in\mathbb{N}$,
\begin{align*}
\bigg[\int_0^\infty |Tx_n(t)|^p~dt\bigg]^{\frac{1}{p}}&=\bigg[\int_{0}^\infty| P_m(t)x_n(t+a)|^p~dt\bigg]^{\frac{1}{p}}\\
          &=\bigg[\int_{n-a}^{n+1-a}|P_m(t)|^p~dt\bigg]^{\frac{1}{p}}\\
          &\geqslant P_m(n+1-a)\\
          &=a_m(n+1-a)^m+a_{m-1}(n+1-a)^{m-1}+\cdots+a_0\\
          &\rightarrow\infty,~n\rightarrow\infty,
\end{align*}
which implies that $T$ is unbounded.

Secondly, we prove that the operator $T^l$ is closed for each $l\in\mathbb{N}$. Let $(x_n)_{n\in\mathbb{N}}$ be such that
$$D(T^l)\ni x_n\rightarrow x\in L_p(0,\infty),~n\rightarrow\infty,$$
and$$T^lx_n\rightarrow y\in L_p(0,\infty),~n\rightarrow\infty.$$
Then, since the function sequences $(x_n(\cdot))_{n\in\mathbb{N}}$ and $(T^lx_n(\cdot))_{n\in\mathbb{N}}$ be convergent in p-norm, then by the Riesz theorem, there exist subsequences $(x_{n_k}\small(\cdot\small))_{k\in\mathbb{N}}$ and $(T^lx_{n_k}(\cdot))_{k\in\mathbb{N}}$ such that $$x_{n_k}(t)\rightarrow x(t),~k\rightarrow\infty~a.e.$$
and $$T^lx_{n_{k}}(t)\rightarrow y(t),~k\rightarrow\infty,~a.e.$$
Thus, $$T^lx_{n_k}(t)=\prod_{i=0}^{l-1}P_m(t+ia)x_{n_k}(t+la)\rightarrow\prod_{i=0}^{l-1}P_m(t+ia)x(t+la),~a.e.$$
By the completeness of Lebesgue measure, we imply that $$\prod_{i=0}^{l-1}P_m(t+ia)x(t+la)=y(t),~a.e.$$
Which implies that $x\in D(T^l)$ and $T^lx=y.$
Therefore, we have already proved that the closed of $T^l$.

Let $\epsilon$ be the set of the equivalence classes represented by p-integrable eventually zero function. Note that $$\epsilon\subseteq C^\infty(T).$$
consider, $$\epsilon\ni x\mapsto Sx\in\epsilon,$$
where the class $Sx$ is represented by $$(Sx)(t):=\begin{cases}\frac{1}{P_m(t-a)}x(t-a),~~~&t>a,\\0,&otherwise\end{cases}$$ then $Sx(\cdot)$ is eventually zero and
\begin{align*}
\int_0^\infty|Sx(t)|^p~dt&=\int_a^\infty\bigg|\frac{1}{P_m(t-a)x(t-a)}\bigg|^p~dt\\
                         &=\int_0^\infty\bigg|\frac{1}{P_m(t)}x(t)\bigg|^p~dt\\
                         &\leqslant\int_{0}^{\infty}\frac{1}{|P_m(0)|^p}|x(t)|^p~dt\\
                         &=\frac{1}{a_0}\int_{0}^{\infty}|x(t)|^p~dt\\
                         &=\frac{1}{a_0}\|x\|_p^p<+\infty.
\end{align*}
Thus the mapping $S:\epsilon\rightarrow\epsilon$ is well defined.

Now, we prove that $TS=I$. Since $\forall x(t)\in\epsilon$, we have that
\begin{align*}
TS(x(t))&=T\bigg(\frac{1}{P_m(t-a)}x(t-a)\bigg)\\
        &=P_m(t)\frac{1}{P_m(t)}x(t)\\
        &=x(t).
\end{align*}
Thus $$TS=I.$$
Let $x\in\epsilon$ represented by $x(\cdot)$ be arbitrary and $supp~x$ be the support of $x(\cdot)$ i.e $\{t\in(0,\infty)|~x(t)\neq0\}$.
Since $x(\cdot)$ is eventually zero, $\exists N\in\mathbb{N}$, $$supp~x\subseteq[0,Na].$$
We can easily show that $\forall x\in \epsilon$, $$T^nx\rightarrow 0(n\rightarrow\infty).$$
Since by the definition of $T,~\forall n\geqslant\mathbb{N}$,
$$T^nx(t)=\prod_{i=0}^{n-1}P_m(t+ia)x(t+na)=0.$$
Next, we will prove that $$\forall x\in \epsilon,~S^nx\rightarrow0~(n\rightarrow\infty).$$
When $\forall x\in\epsilon$ and $\forall n\in\mathbb{N}$, we can conclude that
$$(S^nx)(t)=\begin{cases}\frac{1}{\prod\limits_{i=1}^nP_m(t-ia)}x(t-na),~&t>na,\\0,&otherwise.\end{cases}$$
Since $$\frac{1}{|\prod\limits_{i=1}^nP_m(na-ia)|^p}\rightarrow0,~n\rightarrow\infty,$$
we have that
\begin{align*}
\|S^nx\|_p^p&=\int_0^\infty|S^nx(t)|^p~dt\\
            &=\int_{na}^{\infty}\bigg|\frac{1}{\prod\limits_{i=1}^{n}P_m(na-ia)}x(t-na)\bigg|^p~dt\\
            &\leqslant\frac{1}{|\prod\limits_{i=1}^{n}P_m(na-ia)|^p}\int_{na}^{+\infty}|x(t-na)|^p~dt\\
            &\leqslant\frac{1}{|\prod\limits_{i=1}^{n}P_m(na-ia)|^p}\|x\|_p^p\\
            &\rightarrow0,~n\rightarrow\infty.
\end{align*}
Therefore, by the Sufficient~Condition~for~Hypercyclicity(), the operator $T$ is $hypercyclic$.

Now, we are to show that $T$ has a dense set of periodic point. Since $$C_0(0,\infty)\subset L_p(0,\infty)$$ and $C_0$ is dense in $L_p$,
we just show that $per(T)$ is dense in $C_0$.

Indeed, for any $x\in L_p(0,Na)$ and $\forall N\in\mathbb{N}$, we define the $x_N\in L_p(0,\infty)$ with the following representative
$$x_N(t):=\frac{1}{\prod\limits_{i=1}^{kN}P_m(t-ia)}x(t-kNa),~t\in[kNa,(k+1)Na),~k\in\mathbb{Z_+}.$$
Since
\begin{align*}
T^Nx_N(t)&=\prod\limits_{i=0}^{N-1}P_m(t+ia)\frac{1}{\prod\limits_{i=1}^{kN}P_m(t+Na-ia)}x(t-kNa+Na)\\
         &=\frac{1}{\prod\limits_{i=Na+1}^{kN}P_m(t-ia)}x(t-kNa+Na)\\
         &=\frac{1}{\prod\limits_{i=1}^{(k-1)N}P_m(t-ia)}x(t-(k-1)Na)\\
         &=x_N(t),
\end{align*}
we have that $x_N(t)$ is a periodic point of $T$ with period $N$.

In fact, with $$D_k:=[kNa,(k+1)Na),~k\in\mathbb {Z_+},$$
we get that
\begin{align*}
\int_{0}^{\infty}|x_N(t)|^p~dt&=\sum_{k=0}^{\infty}\int_{D_k}\bigg|\frac{1}{\prod\limits_{i=1}^{kN}P_m(t-ia)}x(t-kNa)\bigg|^p~dt\\
                             &\leqslant\sum_{k=0}^{\infty}\int_{D_k}\bigg|\frac{1}{(a_0)^{kN}}\bigg|^p|x(t-kNa)|^p~dt\\
                             &=\sum_{k=0}^\infty\frac{1}{a_0^{kNp}}\int_{0}^{Na}|x(t)|^p~dt<\infty
\end{align*}
by $a_0>1$ and $x\in L_p(0,Na)$. Thus for any $N\in\mathbb{N}$, $x_N\in L_p(0,\infty).$

Let $y\in C_0(0,\infty)$, then $$\exists n\in\mathbb{N},~y(t)=0,~t>na.$$
Let $y_N$ be the periodic point of $T$ of an arbitrary period $N\geqslant n$ defined based on $y(\cdot)$ as in ( ).

Now, we show that $per(T)$ is dense in $C_0(0,\infty)$. Since by $a_0>1$ and $y\in L_p(0,\infty)$, we have that
\begin{align*}
\|y-y_N\|_P^P&=\int_{0}^{\infty}|y(t)-y_N(t)|^p~dt\\
             &=\sum_{k=1}^{\infty}\int_{D_k}\bigg|\frac{1}{\prod\limits_{i=1}^{kN}P_m(t-ia)}y(t-kNa)\bigg|^p~dt\\
             &\leqslant\sum_{k=1}^{\infty}\frac{1}{a_0^{kNp}}\|y\|_p^p  \rightarrow0,~N\rightarrow\infty.
\end{align*}

Thus which implies that the set $per(T)$ of periodic points of $T$ is dense in $L_p(0,\infty)$, and hence, the operator $T$ is chaotic.
\end{proof}
\section{Proof of Theorem \ref{3}}\label{S5}
In this section, we will prove Theorem \ref{3}.
\begin{proof}
We first show that the operator $T$ is unbounded. Consider $x_n\in D(T)$ such that $\|x_n\|_\infty=1$.

Let the function
$$x_n=\begin{cases}1,~&t\in[0,na),\\ \frac{1}{w^2(t-na)},~&t\geqslant na\end{cases}$$
Since
$$\|Tx_n\|_\infty=\max\limits_{t\geqslant0}|w(t)x_n(t+a)|\geqslant\max_{t\geqslant(n-1)a}\bigg|\frac{w(t)}{w(t-(n-1)a)}\bigg|\rightarrow\infty,~n\rightarrow\infty.$$
Thus, the operator $T$ is unbounded.

Next, we are to show that $T^l$ is closed for $\forall m\in\mathbb{N}$. Let a sequence $(x_n)_{n\in\mathbb{N}}$ be such that
$$D(T^l)\ni x_n\rightarrow x\in C_0[0,\infty),~n\rightarrow\infty,$$
and
$$T^lx_n\rightarrow y\in C_0[0,\infty),~n\rightarrow\infty. $$
Hence, $\forall~t\geqslant0,$,
$$T^lx_n(t)=\prod_{i=0}^{l-1}x_n(t+la)\rightarrow\prod_{i=0}^{l-1}w(t+la)x(t+la),~n\rightarrow\infty.$$
Which implies that $$\prod_{i=0}^{l-1}w(t+ia)x(t+la)=y(t),~t\geqslant0.$$
Therefore $x\in D(T^l)$ and $T^lx=y$. So $T^l$ is closed.

Let $\varepsilon$ be the dense in $C_0[0,\infty)$ of continuous function with compact support.

Consider
$$\epsilon\ni x\mapsto Sx\in\epsilon,$$
where $Sx$ is defined by $$S(x)(t)=\begin{cases}\frac{x(0)}{a}t,&t<a,\\\frac{1}{w(t-a)}x(t-a),&t\geqslant a.\end{cases}$$
The mapping $S:\epsilon\rightarrow\epsilon$ is well defined, and $\forall x\in\epsilon$,
$$TSx(t)=w(t)Sx(t+a)=w(t)\frac{x(t)}{w(t)}=Ix(t).$$
Thus we get that $TS=I$.

Now, we prove that for $\forall x\in\epsilon$,$$T^nx\rightarrow0,~n\rightarrow\infty.$$
Indeed, let $\forall x\in\epsilon$, then there exists $N\in\mathbb{N}$, $$supp~x\subset [0,Na].$$
By the definition of $T$,$$\forall n\geqslant N,~T^nx=0,$$
which implies that $T^n\rightarrow0,~n\rightarrow\infty.$

Next, we will show that $$S^n\rightarrow0,~n\rightarrow\infty.$$
Inductively, for any $n\in\mathbb{N}$,
$$(S^n)(t)=\begin{cases}0,~~&0\leqslant t<(n-1)a,\\ \frac{1}{w(t-a)w(t-2a)\cdots w(t-na)}\frac{x(0)}{a}(t-(n-1)a),~~&(n-1)a\leqslant t<na,\\\frac{1}{w(t-a)w(t-2a)\cdots w(t-na)}x(t-na),~~&t\geqslant na.\end{cases}$$
Since
\begin{align*}
\|S^nx\|_{\infty}&=\sup_{t\geqslant0}|(S^nx)(t)|\\
                 &\leqslant\frac{1}{b^n}x(t-na)\rightarrow0,~n\rightarrow\infty.
\end{align*}
So, by Theorem \ref{chaoxunhuan}, we know $T$ is $Hypercyclicity$.
\end{proof}
We find that the operators given in the \cite{jimenez2022linear} are examples of our thesis.

\bibliographystyle{unsrt}
\bibliography{cite}

\end{document}